\documentclass[11pt]{amsart}
\usepackage{amsmath}
\usepackage{amsxtra}
\usepackage{amscd}
\usepackage{amsthm}
\usepackage{amsfonts}
\usepackage{amssymb}
\usepackage{eucal}
\newtheorem{cor}[subsection]{Corollary}
\newtheorem{lem}[subsection]{Lemma}
\newtheorem{prop}[subsection]{Proposition}

\newtheorem{thm}[subsection]{Theorem}
\newtheorem{defn}[subsection]{Definition}
\theoremstyle{remark}

\newcommand{\nc}{\newcommand}

\nc{\on}{\operatorname}
\nc{\CC}{\mathbb{C}}
\nc{\BZ}{\mathbb{Z}}
\nc{\CV}{\mathcal{V}}
\nc{\CO}{\mathcal{O}}
\nc{\CZ}{\mathcal{Z}}
\nc{\CL}{\mathcal{L}}
\renewcommand{\CD}{\mathcal{D}}
\nc{\CP}{\mathcal{P}}
\nc{\Hom}{\on{Hom}}
\nc{\Sym}{\on{Sym}}
\nc{\Spec}{\on{Spec}}
\nc{\map}{\longrightarrow}
\nc{\til}{\widetilde}

\nc\ssec{\subsection}

\begin{document}
\title{Toric arc schemes and quantum cohomology of toric varieties.} 

\author{Sergey Arkhipov and Mikhail  Kapranov}
\thanks{S. Arkhipov's research is partially supported by an NSF grant}

\address{ Department of Mathematics, Yale University, New Haven CT 06520}
\email  {serguei.arkhipov@yale.edu, mikhail.kapranov@yale.edu}
\maketitle
\vskip .5cm
\section {Introduction}

\vskip .3cm

This paper is a part of a larger project devoted to the study of Floer
cohomology in algebro-geometic context, as a natural cohomology theory
defined on a certain class of ind-schemes. Among these ind-schemes
are algebro-geometric models of the spaces of free loops. 

Let $X$ be a complex projective variety. 
Heuristically, $HQ(X)$, the 
quantum cohomology of $X$, 
is a version of the Floer cohomology of the universal cover of the free loop 
space ${\mathcal L}X = \on{Map}(S^1, X)$. However, 
truly infinite-dimensional approaches to $HQ$ are few. One of the most 
interesting is the unpublished 
result of D. Peterson: for a reductive group $G$
the ring  $HQ(G/B)$ is identified with the usual topological 
cohomology {\it ring} 
of $\mathcal F$, the affine flag variety of $G$. To be precise, if 
$A= H_2(G/B)$ is 
the coweight lattice, 
then $H^\bullet({\mathcal F})$ is naturally an algebra over the
 semigroup ring ${\mathbb  C}[A_+]$ 
of the dominant cone in $A$ and 
$$HQ(G/B) = H^\bullet({\mathcal F}) \otimes_{{\mathbb C}[A_+]}
 {\mathbb C}[A].$$ 
The ind-scheme $\mathcal F$ has the homotopy type
 of ${\mathcal L}G$. Peterson's 
approach was 
deemed mysterious as the analog of $\mathcal F$ for varieties not
 of the form $G/P$ 
was not clear.  

In the present paper we consider the case when $X$ is a toric Fano variety. 
It turns out that the concept of arc schemes (see \cite{DL} \cite{KV}) provides
a replacement of the space $\mathcal F$ in this situation. 
The main  result of the paper, Theorem \ref{main},  identifies
$HQ(X)$ with an appropriate localization of the topological cohomology
of $\Lambda^0X$,  the toric version for the arc scheme.
The action of the generators $q$ arises from the multiplication
of arcs by the parameter $t$. For example, when $X=P^{N-1}$,
the scheme $\Lambda^0X$ is the infinite-dimensional schematic projective
space
$$P({\mathbb C}^N[[t]]) = \on{Proj} \, {\mathbb C}[z_{i,n},\, 
 i=1, ..., N, \,
n\geq 0],$$
where the  $z_{i,n}$ are the coefficients of $N$ indeterminate formal
power series $z_i(t) = \sum_{n=0}^\infty z_{i,n} t^n$.
 Its topological cohomology ring  is ${\mathbb C}[q]$ where $q$
is the hyperplane class. The simultaneous multiplication of the series by $t$
embeds this scheme into itself as a projective subspace of
complex codimension $N$, giving the familiar relation $q=x^N$
in $HQ(P^{N-1})$. 

The scheme $\Lambda^0X$ is a part of an ind-scheme $\Lambda X$
with an action by $H_2(X, \BZ)$ which can be considered as
an algebro-geometric model of the universal cover of $\CL X$. 
For example, for $X=P^{N-1}$
$$\Lambda (P^{N-1}) = P\bigl( \CC ^N((t))\bigr) = 
``\lim \limits _{\longrightarrow} ``  _m \, \on{Proj} \,
\CC[z_{i,n}, i=1, ..., N, \, n\geq -m].$$
Note that a similar model was constructed by Y. Vlassopoulos \cite{Vlas}
(see also \cite{Ir})
who developed the approach of Givental \cite{Giv}. The Vlassopoulos
 model $V(X)$
is a union of finite-dimensional varieties while $\Lambda X$ is not.
For example, for $X=P^{N-1}$ we have 
$$V(P^{N-1}) = P(\CC^N[z, z^{-1}]) = ``\lim\limits_{\longrightarrow}``{} _m
``\lim\limits_{\longrightarrow}`` {} _M \on {Proj} \, \CC[z_{i,n}, i=1, ..., N,
\, \-m\leq n\leq M].$$
Our space $\Lambda X$, constructed as it is from infinite-dimensional
schemes, is locally an ind-pro-object in the category of finite-dimensional
varieties. The duality between the ind- and pro-directions is essential
in our approach to Floer cohomology. For example, in the study of the space
$V(X)$ in {\it loc. cit.} the indices $m$ and $M$ are treated differently
while from the point of view of $\CC[z, z^{-1}]$ they are completely
symmetric. In our approach this difference is built in the
definition of an ind-scheme. 

The idea of modeling the loop spaces by considering formal Laurent series
as opposed to Laurent polynomials was advocated for a long time by
V. Drinfeld. In fact, our ind-scheme $\Lambda X$ is an example of
a T-smooth ind-scheme in the sense of \cite {Dr}. 
We would like to thank H. Iritani for bringing our attention to his work
and to the papers  of Givental and Vlassopoulos.

\vskip .5cm

\section{Quantum cohomology of toric varieties (after Batyrev).}

\vskip .3cm
\subsection{Combinatorics of toric varieties.}
Let $X$ be a smooth compact toric variety defined over $\CC$. From now on
 we fix the realization of $X$ as a quotient, 
following M. Audin and D. Cox, see  \cite{Au},
\cite{Cox}. More precisely, we represent $X=(\CC^N-\CD)/S$
where  $S$ is  a torus acting diagonally on $\CC^N$ (with coordinates $z_1,
 \ldots,z_N$) and  
 $\CD\subset \CC^N$ is a closed subscheme called
the exceptional locus. It has the following structure (see \cite{Bat}
for more detail). Denote the coordinate hyperplane
 $\{ z_i=0\}\subset \CC^N$ by $V_i$.  
For $I\subset\{1,\ldots,N\}$ denote
$
V_I:=\underset{i\in I}{\bigcap}V_i
$. The set $\CD$ has the form 
$$\CD=\underset{I\in\CP}\bigcup V_I$$
where $\CP\subset 2^{\{1,\ldots,N\}}$ is a certain subset
whose elements are called primitive collections. 
It is known that $S$ acts freely on $\CC^N-\CD$ and $X$ is the categorical
quotient of $\CC^N-\CD$ by $S$.

We denote by $A$ (resp. by $B$) the lattice of 1-parametric subgroups
in $S$ (resp.,
the lattice of characters of $S$). In other words,
$$
A=\Hom(\CC^*,S),\ B=\Hom(S,\CC^*).
$$
By construction we have 
$$B\til\to\on{Pic}(X)=H^2(X,\BZ),\ A=B^\vee\til\to H_2(X,\BZ).$$

Below we denote the  image of $V_i-\CD$ in  $X$ by $Z_i$. Notice
that  $Z_i$ is a divisor in $X$. 
Consider the class of $Z_i$ in $ \on{Pic}(X)$ 
denoted by $[Z_i]$. In terms of the isomorphism $\on{Pic}(X)=\Hom(S,\CC^*)$
 the 
class $[Z_i]$ corresponds to the character of $S$ as follows:
$$
S\overset{\on{\tiny{action}}}\map(\CC^*)^N\overset{p_i}\map\CC^*,
$$
where $p_i$ denotes the projection to the $i$-th factor.

\subsection{Combinatorial description for  cohomology of $X$.} 
The following statement
is the combinatorial description of the $\CC$-algebra $H^\bullet(X,\CC)$.

\begin{prop}
The $\CC$-algebra $H^\bullet(X,\CC)$ has generators $[Z_1],\ldots,[Z_N]\in
H^2(X,\CC)$ and the relations as follows:
\begin{itemize}
\item[(i)]
Linear relations among the $[Z_i]$ in $H^2(X,\CC)$.
\item[(ii)]
For any $I\in \CP$ we have $\underset{i\in I}\prod [Z_i]=0$.\qed
\end{itemize}
\end{prop}

\subsection{Combinatorial description for  quantum cohomology of $X$
 (after Batyrev).}
Below we present the main result of Batyrev
(see \cite{Bat}, Theorem 9.5) as made more precise
in \cite{CoxKatz}, Example 11.2.5.2,  providing the combinatorial 
description of the quantum cohomology algebra 
for a toric Fano  variety $X$.

Namely, as a vector space, the quantum cohomology algebra for $X$ denoted by
$HQ(X)$ equals {\em by definition} $H^\bullet(X,\CC)\otimes\CC[A]$.
 Here $\CC[A]$ 
denotes the group algebra of the lattice $A$ formed by Laurent polynomials
$\underset{a\in a}\sum c_a q^a$.

The multiplication on $HQ(X)$ as on a $\CC[A]$-algebra is defined by
 generators and 
relations as follows.

The diagonal action of $S$ on $\CC^N$ defines a homomorphism
$$\beta: A\to \BZ^N = \on{Hom}(\CC^*, (\CC^*)^N).$$
 Let $A_+=\beta^{-1}(\BZ^N_+)$.

\begin{prop} \label{relations}Assume that $X$ is a Fano variety. 
The $\CC[A]$-algebra $HQ(X)$ has generators $[\CZ_1],\ldots,[\CZ_N]$
corresponding to the cocycles $[Z_1],\ldots,[Z_N]
\in H^2(X,\CC)$ and the relations as follows:
\begin{itemize}
\item[(i$)_q$] The generators $[\CZ]_i$ satisfy the same linear
 relations as the 
elements $$[Z_1],\ldots,[Z_N]
\in H^2(X,\CC).$$
\item[(ii$)_q'$] For every $a\in A_+$ we have 
$$\prod_{i=1}^N[\CZ_i]^{\beta_i(a)}=q^a.
\qed
$$
\end{itemize}
\end{prop}

\vskip 0.1pt
\noindent
{\bf Remark:} In particular, it follows
that the $\CC[A]$-algebra defined by the relations (i$)_q$ and 
(ii$)_q'$ is free as a $\CC[A]$-module and has rank equal to
$\on{dim}\ H^\bullet(X,\CC)$.

\vskip .5cm
\section{Cohomology of the toric arcs space.}\label{arcs}

\vskip .3cm

Recall (see  \cite{DL}, \cite{KV} ) that for every $\CC$-scheme $Y$ we
 have the scheme 
$L^0Y$ of formal arcs, representing the following functor:
$$
\Hom_{\mathsf{Sch}}(U,L^0Y)=\Hom_{\mathsf{LRS}}((U,\CO_U[[t]]),(Y,\CO_Y))
$$
where $\mathsf{LRS}$ is the category of locally ringed spaces. For $Y=\Spec(A)$
affine and $U=\Spec(R)$ we have 
$$
\Hom_{\mathsf{Sch}}(U,L^0Y)=\Hom_{\CC\on{-}\mathsf{Alg}}(A,R[[t]]).
$$
It follows that 
$L^0Y = \on { Spec} \, A^{[[t]]}$ where $A^{[[t]]}$ is generated by symbols 
$a[n], a\in A, n\in {\mathbb Z}_+$ 
subject to 
$$(ab)[n] = \sum_{i+j=n} a[i] b[j],$$
compare with the formulas of Borisov \cite{Bor}.  

In particular, for $Y=\CC^N$ with coordinates 
$z_1,\ldots,z_n$ the set $L^0Y(\CC)$ is identified with $\CC^N[[t]]$, the set
 of 
$N$-tuples of formal Taylor series
$$
(z_1(t),\ldots,z_N(t)),\ z_i(t)=\sum_{n=0}^\infty z_{i,n}t^n.
$$ 
We assume the situation and notation of Section 2.
\begin{prop}
The torus $S$ acts freely on the scheme $L^0Y-L^0\CD$.\qed
\end{prop}
\begin{defn}
The Toric arcs scheme $\Lambda^0X$ is defined to be the categorical quotient
$\Lambda^0X:=(L^0Y-L^0\CD)/S$.
\end{defn}

\vskip 0.1pt
\noindent
{\bf Remarks:} (a) The scheme $\Lambda^0X$ is not quasicompact, i.e.,
it cannot be covered by finitely many affine open parts. 
Its construction is completely analogous to that of toric varieties
in Section 2, so it can be considered as an ``infinite-dimensional
toric variety''. In fact, it is acted upon by the group scheme
$$(\CC ^*)^\infty = \on{Spec} \CC[z_{i,n}^{\pm 1}].$$
The subgroup $S$ embedded diagonally into $(\CC ^*)^\infty$,
acts trivially, so the role of the torus acting on a toric
variety, is played in this case by $(\CC ^*)^\infty/S$.

(b) Let $\mathbb{A}^\infty:=\Spec\CC[x_1,x_2,\ldots]$.
Recall \cite {KT} that a scheme $Y$ is called essentially smooth
if it is covered by open subsets which are isomorphic to
$W\times \mathbb{A}^\infty$ where $W$ is a smooth algebraic
variety. In the sequel we will call such schemes
simply smooth. The scheme $\Lambda^0X$ is smooth, in fact
it is covered by open subsets isomorphic to $\mathbb{A}^\infty$. 

\vskip .1cm

Let $p:\ L^0\CC^N\map\CC^N$ be the canonical projection. We define 
$\CV_i=p^{-1}(V_i)\subset L^0\CC^N$ and $\CZ_i\subset\Lambda^0X$ to be the 
image of
$\CV_i-L^0\CD\cap\CV_i$.

Following \cite{KT}, for any $\CC$-scheme (possibly, of infinite type)
$\Sigma$ we introduce the complex topology on the set $\Sigma(\CC)$. In
 particular we
can speak about the topological cohomology
$H^\bullet(\Sigma(\CC),\CC)$. The latter space will be shortly denoted by  
$H^\bullet(\Sigma,\CC)$.

\begin{prop}
Odd cohomology spaces of $\Lambda^0X$ vanish.  The graded ring 
$H^{2\bullet}(\Lambda^0X,\CC)$ is isomorphic to the symmetric algebra 
$\Sym^\bullet(H^2(X,\CC))$.
\end{prop}
\begin{proof}
We begin the proof with the following statement.

\begin{lem}
$H^\bullet(L^0\CC^N-L^0\CD,\CC)$ vanishes in positive degrees.
\end{lem}

\begin{proof}
Since $\CD=\underset{I\in\CP}\bigcup V_I$, we have $L^0\CD=
\underset{I\in\CP}\bigcup L^0V_I$ and each $L^0V_I$ is an affine subspace in 
$L^0\CC^N$ of infinite codimension. The statement of the Lemma follows immediately.
\end{proof}

Now, in view of the previous Lemma, the Proposition follows from the 
Serre spectral sequence of the fibration
$\rho:\ L^0\CC^N-L^0\CD\overset{S}\map\Lambda^0X
$
and the fact that $H^\bullet(S,\CC)$ is the exterior algebra of the space 
$B\otimes \CC=H^2(X,\CC)$.
\end{proof}

\begin{cor}
The classes $[\CZ_i]\in H^2(\Lambda^0X,\CC)$ generate $ H^\bullet
(\Lambda^0X,\CC)$
as a ring.\qed
\end{cor}

\vskip .5cm

\section{Self-embeddings and Gysin maps.}

\vskip .3cm
Let $a\in A_+$. Then the 1-parameter subgroup $t^{\beta(a)}$ of $S$ can
 be considered 
as a $\CC$-point of $L^0(\CC^*)^N$. Using the action of $(\CC^*)^N$ on $\CC^N$,
for any $\CC$-point $\gamma(t)$ of $L^0\CC^N$ 
we can form a new $\CC$-point $t^{\beta(a)}\gamma(t)$. It is clear that the 
correspondence $\gamma(t)\mapsto t^{\beta(a)}\gamma(t)$ gives rise to a 
morphism of schemes
$\epsilon_a:\ L^0\CC^N\map L^0\CC^N$ and that $\epsilon_a$ is a closed embedding. 
It is also clear (from the fact that $\beta$ is a homomorphism) that
$$
\epsilon_{a+b}=\epsilon_{a}\circ\epsilon_{b},\ a,b\in A_+.
$$
We define the  subscheme $L^a\CC^N\subset L^0\CC^N$ to be the image of 
$\epsilon_{a}$.

\begin{prop} \label{intersection}
$L^a\CC^N$ is equal to the closed subscheme  of $L^0\CC^N$ formed by 
$N$-tuples of series $\gamma(t)=(\gamma_1(t),\ldots,\gamma_N(t))$
such that $\on{ord}_t\gamma_i(t)\ge\beta_i(a)$.\qed
\end{prop}

Introduce a partial order $\le$ on $A_+$ by putting $a\le b$ iff $b-a\in A_+$. 
It is clear from the above Proposition
 that $L^b\CC^n\subset L^a\CC^N$ if and only if $a\le b$. 
\begin{cor}
We have $\on{codim}_{L^a\CC^N}L^b\CC^N =\sum_{i=1}^N\beta_i(b-a)$.\qed
\end{cor}

%We define a
%locally closed subscheme
%$L{=a}\CC^N\subset L^0\CC^N$ by
%$$
%L^{=a}\CC^N:=L^a\CC^N-\underset{b>a}\bigcup L^b\CC^N.
%$$
It is clear that for $a\in A_+$ we have 
$\epsilon_{a}(L^0\CC^N-L^0\CD)\subset L^0\CC^N- L^0\CD$ and that $\epsilon_{a}$ 
commutes with the action of $S$. Therefore the map descends to a closed
embedding
$$
E_a:\ \Lambda^0X\map \Lambda^0 X, \ E_{a+b}=E_a\circ E_b.
$$
The image of $E_a$ will be denoted by $\Lambda^aX$. We also conclude that 
$\Lambda^bX\subset\Lambda^aX$ if and only if $a\le b$ and, if so, 
$$
\on{codim}_{\Lambda^aX}\Lambda^bX =\sum_{i=1}^N\beta_i(b-a).
$$
%We also define $\Lambda^{=a}X:=\Lambda^aX-\underset{b>a}\bigcup \Lambda^bX$.

\begin{prop} \label{id}
The inverse image map $E^*_a:\ H^\bullet(\Lambda^0X,\CC)\to
H^\bullet(\Lambda^0X,\CC)$
is the identity for any $a\in A_+$.
\end{prop}

\begin{proof}
The maps $\epsilon_a$ and $E_a$ define an embedding of the fibration
$$\rho:\ L^0\CC^N-L^0\CD\overset{S}\map\Lambda^0X$$ into itself. On the 
cohomology of the total spaces $E^*_a=\on{Id}$. On the cohomology of fibers 
$E^*_a$ is identity as well. Our statement follows from functoriality of
Serre spectral sequences.
\end{proof}

Let $Y\overset{q} \to Z$ be a smooth closed embedding of (possibly
infinite-dimensional)  schemes of
finite codimension $d$. This implies that locally on the complex
topology the embedding of $\CC$-points is identified with
$Y(\CC) \to Y(\CC)\times \CC^d$. This gives an identification in the
derived category
$$\underline{R\Gamma}_{Y(\CC)} (\underline{\CC}_{Z(\CC)}) \to
 \underline{\CC}_{Y(\CC)}
[2d]$$ and hence gives the Gysin maps
$$q_!: H^m(Y, \CC) \to H^{m+2d}(Z,\CC).$$

\vskip .2cm

Let $\CC[A_+]$ be the semigroup algebra of $A_+$. Its elements are 
represented as (formal) polynomials $\underset{a\in A_+}\sum c_aq^a$ in a 
variable
$q\in\Hom(A,\CC^*)$.
We make the space $H^\bullet(\Lambda^0X,\CC) $ into a 
$\CC[A_+]$-module by putting for $\alpha\in H^\bullet(\Lambda^0X,\CC)$
\begin{equation} \label{formula}
\left(\sum c_aq^a\right)\cdot\alpha:=\sum c_a E_{a\hspace{1pt}!}(\alpha).
\end{equation}
Here $ E_{a\hspace{1pt}!}:\ H^m(\Lambda^0X,\CC)\map 
H^{m+2\sum\beta_i(a)}(\Lambda^0X,\CC)$
is the Gysin map corresponding to the smooth closed embedding $E_a$. The fact that 
the formula above defines a $\CC[A_+]$-module structure follows immediately from the 
composition rule $E_{a+b}=E_a\circ E_b$ and the next Lemma.

\begin{lem}
Let $Y_1\overset{q}\to Y_2\overset{p}\to Y_3$ be smooth closed embeddings
of schemes of finite codimension.  Then we have the equality of Gysin maps
$$
(pq)_!=p_!\circ q_!.\qed
$$
\end{lem}

Recall the  key property of Gysin maps (projection formula). It 
follows at once from the sheaf-theoretic definition of Gysin maps above.

\begin{prop}
Let 
$q:\ Y\map Z$ be a closed embedding of smooth $\CC$-schemes of 
finite codimension $d$. Then the Gysin map
$$
q_!:\ H^\bullet(Y, \CC)\map H^\bullet(Z, \CC)
$$
is $H^\bullet(Z, \CC)$-linear, i.e. we have 
$$
q_!(q^*(\alpha)\cap\beta)=\alpha\cap q_!(\beta), \
\alpha\in H^\bullet(Z, \CC),\ \beta\in H^\bullet(Y, \CC).
\qed
$$
\end{prop}

Applying this to $q=E_a$, $Y=Z=\Lambda^0X$ and using Proposition~\ref{id}
we get the following

\begin{cor}
The formula \ref{formula} makes the ring $H^\bullet(\Lambda^0X, \CC)$ into
a $\CC[A_+]$-algebra.\qed
\end{cor}

Now we come to the main statement of the paper.

\begin{thm} \label{main}
The localized algebra
$$
H^\bullet(\Lambda^0X, \CC)\otimes_{\CC[A_+]}\CC[A]
$$
is isomorphic the the quantum cohomology algebra $HQ(X,\CC)$
(isomorphism of $\CC[A]$-algebras).
\end{thm}

\vskip .5cm

\section{Proof of Theorem~\ref{main}.}

\vskip .3cm
We construct a homomorphism of $\CC[A]$-algebras 
$$
\Phi:\ HQ(X)\map H^\bullet(\Lambda^0X, \CC)\otimes_{\CC[A_+]}\CC[A]
$$
by putting $\Phi([Z_i])=[\CZ_i]\otimes 1$. In order to see that this indeed defines 
a homomorphism, we need to check the relations $(i)_q$ and $(ii)'_q$ from 
Proposition~\ref{relations}. For $(i)_q$ this follows form the identification
$$
H^2(X,\CC)\til\to H^2(\Lambda^0X, \CC)\til\to B\otimes\CC
$$
which takes $[Z_i]$ to $[\CZ_i]$.
It remains to check  $(ii)'_q$. 

 As in Section 3, we introduce coordinates $z_{i,n}$, $i=1,\ldots,N$ and $n\in\BZ_+$, in $L^0\CC^N$ 
so that a $\CC$-point has the form  
$$
\gamma(t)=(\gamma_1(t),\ldots,\gamma_N(t)),\
\gamma_i(t)=\sum_{n=0}^\infty z_{i,n} t^n.
$$
Let $\CV_{i,n}\subset L^0\CC^N$ be the hypersurface $z_{i,n}=0$ and 
$\CZ_{i,n}\subset \Lambda^0X$ be the corresponding divisor. In particular we have
$\CZ_i=\CZ_{i,0}$. On the other hand, $\CZ_{i,n}$ is linearly equivalent to 
$\CZ_{i,m}$ for any $n$, $m$. This is because the function $z_{i,n}/z_{i,m}$ on
$L^0\CC^N$ is homogeneous of degree $0$ with respect to the action of the torus $S$
and it descends a rational function on $\Lambda^0X$ with divisor
$\CZ_{i,n}-\CZ_{i,m}$. Therefore we can write the self-intersection
$[\CZ_i]^r$ as the class of the subvariety 
$\CZ_{i,0}\cap\ldots\cap\CZ_{i,r-1}$.
So the class $\prod[\CZ_i]^{\beta_i(a)}$ in $(ii)'_q$ is the class of the subvariety
$$\bigcap_{i=1}^N\CZ_{i,0}\cap\ldots\cap\CZ_{i,\beta_i(a)-1}$$
which coincides with $\Lambda^aX$ in virtue of 
Proposition~\ref{intersection}.

Further, $[\Lambda^aX]$ is equal to the result of action of $q^a=E_{a!}$ on the 
element $1\in H^\bullet(\Lambda^0X,\CC)$. So the homomorphism $\Phi$ is well defined.
Notice that it is surjective since  $H^\bullet(\Lambda^0X,\CC)$ is generated by the 
classes $[\CZ_i]$ as a $\CC$-algebra. It remains to prove that  $\Phi$ is injective.

\begin{prop} \label{free}
The cohomology ring $H^\bullet(\Lambda^0X,\CC)$ is torsion free as a 
$\CC[A_+]$-module and has generic rank equal to $\dim \, H^\bullet(X, \CC)$.
\end{prop}
\begin{proof}
We define a
locally closed subscheme
$L^{=a}\CC^N\subset L^0\CC^N$ by
$$
L^{=a}\CC^N:=L^a\CC^N-\underset{b>a}\bigcup L^b\CC^N.
$$
We also define $\Lambda^{=a}X:=\Lambda^{a}X-\underset{b>a}\bigcup\Lambda^{b}X$.  
For $a\in A_+$ let $I^a:\ \Lambda^{=a}X\hookrightarrow \Lambda^0X$ be the embedding.
Recall that   $p$ denotes  the canonical projection $L^0\CC^N\to\CC^N$.
\begin{lem}
We have the equalities 
$$L^{=0}\CC^N=L^0\CC^N-p^{-1}(\CD)\text{ and }
\Lambda^{=0}X=(L^0\CC^N-p^{-1}(\CD))/S.\qed$$
\end{lem}

In particular, $p:\ L^0\CC^N\to\CC^N$ induces a morphism of schemes
$
P^0:\ \Lambda^0X\map X$.

\begin{lem} \label{lemma2}
$P^0$ is a (Zariski) locally trivial fibration with fiber
 $\mathbb{A}^\infty$. 
In particular, the map
$$(P^0)^*:\ H^\bullet(X,\CC)\til\map H^\bullet(\Lambda^{=0}X,\CC)
$$
is an isomorphism of $\CC$-algebras.\qed
\end{lem}

\begin{cor}
The map
$$(E_a)^*\circ (P^0)^*:\ H^\bullet(X,\CC)\til\map
 H^\bullet(\Lambda^{=a}X,\CC)
$$
is an isomorphism of $\CC$-algebras.\qed
\end{cor}

Consider the Cousin spectral sequence (see \cite{H}, p. 227) corresponding to
the filtration of $\Lambda^0X$ by the $\Lambda^aX$ and the constant
sheaf $\underline{\CC}$ on $\Lambda^0X$. It has
$$E_1 = \bigoplus_{a\in A_+} H^\bullet_{\Lambda^{=a}X}(\Lambda^0X, \CC), 
\quad E_\infty = H^\bullet(\Lambda^0X, \CC).$$
Now, since $\Lambda^aX \subset \Lambda^0X$ is a smooth embedding, we have
$$H^p_{\Lambda^{=a}X}(\Lambda^0X, \CC) = 
H^{p+2\sum \beta(a_i)}(\Lambda^{=a}X, \CC) = 
H^{p+2\sum \beta(a_i)}(X, \CC).$$
In particular, the $E_1$-term is concentrated in even degrees. It follows
that the spectral sequence
 degenerates and provides a certain filtration on 
$H^\bullet(\Lambda^0X,\CC)$ with the key property:
$$
\on{gr} (H^\bullet(\Lambda^0X,\CC))=\CC[A_+]\otimes H^\bullet(X,\CC).
$$
as a $\CC[A_+]$-module. Proposition~\ref{free} is proved.
\end{proof}

We conclude that the $\Phi:\ HQ(X)\to 
H^\bullet(\Lambda^0X,\CC)\otimes_{\CC[A_+]}\CC[A]$
is a map from a {\em free} $\CC[A]$-module to a module  of the same 
generic rank (equal to 
$\on{dim} (H^\bullet(X,\CC))$). Moreover, we have already proved that $\Phi$ 
is surjective. 

It follows that
$ H^\bullet(\Lambda^0X,\CC)\otimes_{\CC[A_+]}\CC[A]$ is a free module.
Indeed, freeness of a finitely generated module
is equivalent to the fact that the rank of its
specialization at any point is the same. In any case, the rank of the
specialization is greater or equal than the generic rank.

Now, since the specialization functor is right exact and
$\Phi$ is surjective, we get that the rank of the specialization
is less or equal to the generic rank thus implying the freeness.  
Now, a surjective morphism between free modules of the same
rank should be injective and so an isomorphism.  

Theorem~\ref{main} is proved.

\vskip .5cm

\section{Concluding remarks} 

\vskip .3cm

\subsection{Relation to Floer cohomology.}
Let us explain the relation of the above construction with
our algebro-geometric approach to Floer cohomology. By inverting the
self-embeddings $E_a, a\in A_+$, we include the scheme $\Lambda^0X$
into an ind-scheme 
$$\Lambda X = \text{''}{\lim\limits_{\underset{a\in A}\longrightarrow}}
\text{''} 
\Lambda^a X$$
on which the group $A= H_2(X, \BZ)$ acts by automorphisms. This is a toric
algebro-geometric model of the universal cover of the free loop
space ${\mathcal L}X$. The Floer cohomology of $\Lambda X$ is defined to be
$$HF(\Lambda X) = {\lim\limits_{\underset{a\in A}\longrightarrow}} \,  
H^{\bullet+ 2d(a)}
(\Lambda^a X, \CC), \quad d(a) =\sum \beta(a_i),$$
where the limit is taken with respect to the Gysin maps. This limit
is identified with the extension of scalars from $\CC[A_+]$ to
$\CC[A]$  as in Theorem \ref{main}.

\subsection{More general varieties.} 
In a forthcoming paper we plan
 to generalize
the above  approach to a large class of smooth projective varieties.
Let $X$ be a smooth projective variety over $\CC$ together with 
a choice of a projective embedding, i.e., represented as $X=\on{Proj}(A)$.
We define $\Lambda^0X = \on{Proj} (A^{[[t]]}$ where the algebra
$A^{[[t]]}$ is defined as in  Section \ref{arcs} and has grading
$\on{deg}(a[n])=\on{deg}(a)$. In other words, if $Y=\on{Spec}(A)$
is the cone over $X$, then $\Lambda^0 X$ is
the categorical quotient of $L^0(Y)-\{0\}$ by $\CC^*$. The multiplication
by $t$ defines a self-embedding $\epsilon: L^0Y\to L^0Y$
which descends to an embedding $E: \Lambda^0X\to\Lambda^0X$. 
The ind-scheme
$$\Lambda X = `` \lim\limits_{\longrightarrow} ``
\biggl\{ \Lambda^0X \overset{E}\map \Lambda^0 X \overset{E}\map...\biggr\}$$
is an algebro-geometric analog of the $\BZ$-cover of $\mathcal{L}X$
corresponding to  $c_1(\mathcal{O}_X(1))\in H^2(X, \BZ)$. 
The case of the universal cover can be treated in a similar way.

While  $\Lambda^0X$ is no longer smooth , it
is nevertheless possible to define an analog of the Gysin maps
and perform the limit construction. This is based on the
Cartesian square
$$\begin{matrix} L^0Y  &\buildrel \epsilon \over\longrightarrow &L^0 Y&\cr 
\big\downarrow&&\big\downarrow&\cr 
\{0\} & \longrightarrow& Y&\end{matrix} $$ 
If $d=\dim(X)$, we have a class $\xi\in H^{2d+2}_{\{0\}}(Y, { 
\CC})$ coming from the fundamental 
class in $H^{2d}(X, {\CC})$.
 From the above square we get a class $\eta\in 
H^{2d+2}_{\epsilon (L^0Y)}(L^0Y, {\CC})$. 
The embedding $E$ has  
codimension $2d+2$ 
and $\eta$, being ${\CC}^*$-equivariant, descends to 
$$\zeta\in H^{2d+2}_{E(\Lambda^0X)}
 (\Lambda^0X, {\CC}) = {\rm Hom} (E^! 
\underline {\CC}, 
\underline {\CC}[2d+2]).$$ 
The map $\zeta$ provides a Gysin map in our non-smooth case. 
We define the partial Floer cohomology (corresponding to
the given projective embedding) as 
$$HF(\Lambda X) = \lim\limits _{\map}  \biggl\{ H^\bullet (\Lambda^0X, \CC)
\overset{E_!}\map H^{\bullet+2d+2}(\Lambda^0X, \CC)
\overset{E_!}\map H^{\bullet+4d+4}(\Lambda^0X, \CC)\overset{E_!}\map ...
\biggr\}.$$
In general, we do not expect a direct relation between $HQ(X)$
and the ring structure on $H^\bullet(\Lambda^0X)$.

Note that the homotopy limit of the diagram 
of complexes of  sheaves on $\Lambda^0X(\CC)$
$$\underline {{\CC}}\to E^! \underline{\CC}[2d+2]\to (E^2)^! 
\underline {\CC} [4d+4]\to ... $$ 
can be seen as a topological dualizing sheaf, cf. the discussion
in \cite{Dr} in a more restricted situation.

\end{document}